\documentclass[12pt]{article}
\usepackage{amssymb}
\usepackage{array}
\usepackage{latexsym}
\usepackage{enumerate}
\usepackage{amsmath}
\usepackage{amsfonts}
\usepackage{amsthm}
\usepackage{graphicx}
\usepackage{psfrag}
\usepackage{comment}
\usepackage[english]{babel}
\usepackage[T1]{fontenc}
\ifx\smallsetminus\undefined
\def\smallsetminus{\setminus}
\fi
\title{Algebraic curves and maximal arcs}
\author{A. Aguglia ${}^*$ \and L. Giuzzi ${}^*$
  \and G. Korchm\'aros
  \thanks{Research supported by  the Italian
    Ministry MURST, Strutture geometriche, combinatoria e loro
    applicazioni.}}

\date{}

\theoremstyle{plain}
\newtheorem{prop}{Proposition}[section]
\newtheorem{theorem}[prop]{Theorem}
\newtheorem{corollary}[prop]{Corollary}
\newtheorem{lemma}[prop]{Lemma}
\theoremstyle{definition}
\newtheorem{remark}[prop]{Remark}

\def\cC{\mathcal C}
\newcommand{\cD}{\mathcal D}

\def\cK{\mathcal K}

\def\cX{\mathcal X}
\def\cK{\mathcal K}

\newcommand{\fF}{\mathfrak{F}}
\newcommand{\PG}{\mathrm{PG}}
\newcommand{\GF}{\mathrm{GF}}
\newcommand{\tr}{\mathfrak{T}}
\newcommand{\oGF}[1]{\overline{\GF(#1)}}
\begin{document}

\maketitle
\begin{abstract}
A lower bound on the minimum degree of the plane algebraic curves
containing every point in a large point--set $\cK$ of the
Desarguesian plane $\PG(2,q)$ is obtained. The case where $\cK$ is
a maximal $(k,n)$--arc is considered in greater depth.


\end{abstract}

\section{Introduction}
In finite geometry, plane algebraic curves of minimum degree
containing a given large point--set $\cK$ in $\PG(2,q)$ have been
a useful tool to investigate combinatorial properties of $\cK$.

When $\cK$ is the whole point--set of $\PG(2,q)$, a trivial lower
bound on the degree of such a plane algebraic curve is $q+1$.
G.~Tallini pointed out that this is attained only when the curve
splits into $q+1$  distinct lines of $\PG(2,q)$, all
passing through the same point. He also gave a complete
classification of the absolutely irreducible curves of degree
$q+2$ containing all points of $\PG(2,q)$; see \cite{tal1,tal2}
and also \cite{abko}. If $\cK$ is the complementary set of a line
in $\PG(2,q)$, then the bound is $q$; see \cite{goppa}.

When $\cK$ consists of all internal points to a conic $\cC$ in
$\PG(2,q)$ with $q$ odd, the above lower bound is $q-1$. The
analogous bound for the set of the external points to $\cC$ is $q$.
These bounds were the main ingredients for recent combinatorial
characterisations of point--sets blocking all external lines to
$\cC$; see \cite{ak,giulietti}.

When $\cK$ is a classical unital of $\PG(2,q)$, with $q$ square,
the minimum degree $d$ of an absolutely irreducible curve $\cC$
through $\cK$ is $d=\sqrt{q}+1$. For non--classical unitals, the
best known bound is $d>2\sqrt{q}-4$; see \cite{HK}.

Our purpose is to find similar bounds for slightly smaller, but still
quite large, point--sets $\cK$, say $|\cK|=qt+\alpha$ with
$t<c\sqrt[3]{q}$ and $0\leq \alpha <q$, where $c$ is a suitable
constant. Since no combinatorial
condition on the configuration of $\cK$ is assumed, we are relying
on techniques and results from algebraic geometry rather than on
the constructive methods  used in the papers cited above.

The main result is that if $q>8t^3-12^2+4t-2\alpha+2$, any plane
algebraic curve $\Gamma$ containing every point of $\cK$
has degree $d\geq 2t$. The hypothesis on the
magnitude of $t$ can be relaxed to $q>16t^2-24t-2\alpha+8$
whenever $q$ is a prime.

In some cases, the bound $d\geq 2t$ is sharp, as the following
example shows. Let $\cK$ be the union of $t$  disjoint
ovals. If $q$ odd, or $q$ is even and the ovals are classical,
then $|\cK|=qt+t$ and $d=2t$. The latter case is known to occur
when $\cK$ is a Denniston maximal arc \cite{D} (or one of the maximal arcs
constructed by Mathon and others, \cite{T74,T80,P,M,HM})
minus the common nucleus of the
ovals.

On the other hand, some refinement of the bound is also possible.
Let $\cK$ be any maximal arc of size $|\cK|=qt+t+1$, that is, a
$(qt+t+1,t+1)$--arc. Theorem \ref{teo1.5} shows that if $q$ is
large enough comparing to $t$, then no plane algebraic curve of
degree $2t$ passes through every point of $\cK$. Therefore, the
minimum degree is at least $2t+1$ for such $t$ and this bound is
attained when $\cK$ is one of the above maximal arcs.



The case $n=4$ is considered in more detail. For $q>2^6$, the
minimum degree is $7$ and this is only attained when $\cK$ is a
Denniston arc and the curve $\Gamma$ of minimum degree splits into
three distinct conics with the same nucleus $N\in \cK$, together
with a line through $N$.

\section{Some background on plane algebraic
  curves over a finite field}
\label{background}


A plane projective algebraic curve $\Delta$ is defined over
$\GF(q)$, but viewed as a curve over the algebraic closure
$\overline{\GF(q)}$ of $\GF(q)$, if it has an affine equation
$f(X,Y)=0$, where $f(X,Y)\in\GF(q)[X,Y]$. The curve $\Delta$ is
absolutely irreducible if it is irreducible over the algebraic
closure $\oGF{q}$. Denote by $N_q$ be the number of non--singular
points lying in $\PG(2,q)$ of an absolutely irreducible plane
curve $\Delta$ of degree $d$. From the Hasse-Weil bound
\begin{equation}
\label{hasse}
N_q\leq q+1 +(d-1)(d-2)\sqrt{q}.
\end{equation}
This holds true when singular points of $\Delta$
lying in $\PG(2,q)$ are also counted; see \cite{LeYe}.

The St\"ohr-Voloch bound depends not only on the degree $d$, but also
on a positive integer, the Frobenius order $\nu$ of $\Delta$;
see \cite{SV}.
This number $\nu$ is either $1$ or $\varepsilon_2$, where
$\varepsilon_2$ is the intersection number $I(P,\Delta\cap \ell)$ of
$\Delta$ with the tangent line $\ell$ at a general point $P\in\Delta$. It
turns out that $\varepsilon_2$ is either $2$, or a power, say $p^h$,
of the characteristic $p$ of the plane, and is the minimum of
$I(Q,\Delta\cap r)$, where $Q$ ranges over the non--singular points of
$\Delta$ and $r$ is the tangent to $\Delta$ at $Q$. If $q=p$, then
$\varepsilon_2=2$, and either $\nu=1$, or $\nu=2$ and $p=2$. With
this notation, the St\"ohr-Voloch bound applied to $\Delta$ is
    \begin{equation}
    \label{svplane}
2N_q\leq \nu(d-3)d+d(q+2).
    \end{equation}
The following algebraic machinery can be used to compute $\nu$.
Let $P=(a,b)$ be a non--singular point of $\Delta$ such that the
tangent line to $\Delta$ at $P$ is not the vertical line through
$P$. The unique branch (or place) centred at $P$ has a local
parametrisation, also called a primitive branch representation,
$$x=a+t,\,y=b+\varphi(t)$$ where $f(x,y)=0$ and
$\varphi(t)=b_kt^k+\ldots$ with $k\geq 1$ is a formal power series
with coefficients in $\overline{\GF(q)}[[t]]$; see
\cite{Seidenberg}. Then, $\nu$ is defined to be the smallest
integer such that the determinant
$$\left|
\begin{array}{cc}
    x-x^q & y-y^q \\
     1 & D_t^{(\nu)}(y)
\end{array}
\right|= \left| \begin{array}{cc}
    a-a^q+t-t^q & b-b^q+\varphi(t)-\varphi(t)^q \\
     1 & D_t^{(\nu)}(\varphi(t))
\end{array}
\right|
$$
does not vanish. Here $D_t$ denotes the $\nu$--th Hasse derivative,
that is,
$$D_t^{(\nu)}(\varphi(t))=\binom{k}{\nu} b_kt^{k-\nu}+\ldots.$$

The above idea still works if osculating conics are used in place
of tangent lines, and, in some cases, the resulting bound improves
\eqref{svplane}. Before stating the result, which is the
St\"ohr-Voloch bound for conics, a further concept from algebraic
geometry is needed. Recall that the order sequence of $\Delta$
with respect to the linear system $\Sigma_2$ of the conics of the
plane is the increasing sequence
$0,\epsilon_1=1,\epsilon_2=2,\epsilon_3,\epsilon_4,\epsilon_5$ of
all intersection numbers $I(P,\Delta\cap \cC)$ of $\Delta$ with
conics at a general point $P$.
The Frobenius $\Sigma_2$--order sequence is the subsequence $\nu_0
= 0, \nu_1, \nu_2,\nu_3,\nu_4$ extracted increasingly from the
$\Sigma_2$--order sequence of $\Delta$, for which the following
determinant does not vanish:
\[\left| \begin{array}{ccccc}
    x-x^q& x^2-x^{2q}&y-y^q& xy-x^qy^q& y^2-y^{2q} \\
    1 & 2x& D_t^{(\nu_1)}(y) & D_t^{(\nu_1)}(xy) &  D_t^{(\nu_1)}(y^2)\\
    0&1 & D_t^{(\nu_2)}(y)& D_t^{(\nu_2)}(xy)& D_t^{(\nu_2)}(y^2)\\
    0 & 0& D_t^{(\nu_3)}(y)& D_t^{(\nu_3)}(xy)&D_t^{(\nu_3)}(y^2)\\
0 & 0& D_t^{(\nu_4)}(y)& D_t^{(\nu_4)}(xy)&D_t^{(\nu_4)}(y^2)
\end{array}
\right|. \]
 Assume that $\deg\Delta\geq 3$. The St\"ohr-Voloch
bound for conics, that is for $\Sigma_2$, is
\begin{equation}
\label{svconic}
  5N_q\leq [(\nu_1+\ldots+\nu_4)(d-3)d+2d(q+5)].
\end{equation}
For more on the St\"ohr-Voloch bound see \cite{SV}.

\section{Plane algebraic curves of minimum degree through all the points
of a given point--set} \label{remarkable}

In this section, $\cK$  stands for a set of $qt+\alpha$ points in
$\PG(2,q)$, with $0\leq t\leq q$ and $0\leq \alpha<q$. Let
$\Gamma$ denote a plane algebraic curve of degree $d$ containing
every point $\cK$.
As already mentioned,  $q+1$ is the minimum degree of a plane
algebraic curve containing every point of $\PG(2,q)$.
Thus, since we are looking for lower bounds on $d$, we will only
be concerned with the case where $d\leq q$.

A straightforward counting argument gives the following result.
\begin{lemma}
\label{lemma11} If $d\leq q$, then $d\geq t$.
\end{lemma}
\begin{proof}
  Since $d\leq q$, the linear components of $\Gamma$ do not contain
all the points of $\PG(2,q)$. Choose a point $P\in\PG(2,q)$
not in any of these linear components.
Each of the $q+1$ lines through $P$ meets $\Gamma$ at most
$d$ distinct points. Thus, $(q+1)d\geq |\cK|$, that is
$$d\geq \frac{qt+\alpha}{q+1}\geq t -\frac{t}{q+1}.$$
Since $t<q+1$, the assertion follows.
\end{proof}
Our aim is to improve Lemma \ref{lemma11}. Write $\fF$ for the set
of all lines of $\PG(2,q)$ meeting $\cK$ in at least $1$ point.
Set $m_0=\min\{|\ell\cap\cK|:\ell\in\fF\}$ and
$M_0=\max\{|\ell\cap\cK|:\ell\in\fF\}$.

\begin{theorem}
\label{teo:main}
Let $\Gamma$ be an algebraic plane curve over
$\oGF{q}$ of minimal degree $d$
which passes through all the points of $\cK$.
If
\begin{equation}
\label{eq:co1}
 q>8t^3-16t^2+2t+4-2m_0(2t^2-5t+2)+2M_0(2t-1),
\end{equation}
then
$\deg\Gamma\geq 2t$. For prime $q$, Condition \eqref{eq:co1} may
be relaxed to
\begin{equation}
\label{eq2}
q>8t^2-16t+8-2\alpha+2M_0(2t-1).
\end{equation}
\end{theorem}
\begin{proof}
We prove that if $d\leq 2t-1$, then \eqref{eq:co1}
does not hold. For $q$ prime we show that also  \eqref{eq2}
is not satisfied. Since $\Gamma$ is not necessarily
irreducible, the following setup is required.
\par
 The curves
$\Delta_1,\ldots,\Delta_l$ are the absolutely irreducible
non--linear components of $\Gamma$ defined over $\GF(q)$,
respectively of degree $d_i$;
$r_1,\ldots,r_k$ are the linear components of $\Gamma$ over
$GF(q)$; $\Xi_1,\ldots,\Xi_s$ are the components of $\Gamma$
which are irreducible over $\GF(q)$ but not over $\oGF{q}$.

The idea is to estimate the number of points in $\PG(2,q)$ that
each of the above components can have.

Let $N_i$ be the number of non--singular points of $\Delta_i$
lying in $\PG(2,q)$. Then, \eqref{svplane} holds for any
$\Delta_i$. Let $\nu^{(i)}$ denote the Frobenius order of
$\Delta_i$. If $\nu= \max \{\nu^{(i)}\mid 1\leq i \leq l\}$ and
$\delta= \sum_{i=1}^l d_i$, then
\begin{equation}
\label{eq4} 2\sum_{i=1}^l N_i \leq \sum_{i=1}^l \nu^{(i)}
d_i(d_i-3)+(q+2)d_i \leq \nu \delta(\delta-3)+(q+2)\delta.
\end{equation}
For $q$ prime, $\nu=1$; see \cite{SV}. Since $\nu^{(i)}=1$ can
fail for $q>p$, an upper bound on $\nu^{(i)}$ depending on $d_i$
is needed. As $\nu^{(i)}\leq \varepsilon_2^{(i)}$, a bound on
$\varepsilon_2^{(i)}$ suffices. Since $N_i>0$ may be assumed,
$\Delta_i$ has a non--singular point $P$ lying in $\PG(2,q)$. If
$\ell$ is the tangent to $\Delta_i$ at $P$, then
\[
d_i=\sum_{Q\in\ell\cap\Delta_i}I(Q,\ell\cap\Delta_i)=%
I(P,\ell\cap\Delta_i)+\sum_{\begin{subarray}{c}Q\in\ell\cap\Delta_i \\
   Q\neq P
 \end{subarray}}I(Q,\ell\cap\Delta_i)\geq
 \varepsilon_2^{(i)}+m_0-1,
\]
whence $\nu^{(i)}\leq d_i-m_0+1$. From \eqref{eq4},
\begin{equation}
\label{eq6}
 2\sum_{i=1}^l N_i\leq \begin{cases}
\delta(\delta-3)+(q+2)\delta&\mbox{$q\geq 3$ prime.}\\
   (\delta-m_0+1)\delta(\delta-3)+(q+2)\delta&\mbox{otherwise.\/}

\end{cases}
\end{equation}

If $\Delta_i$ has $M_i$ singular points, from Pl\"ucker's theorem
$M_i\leq \frac{1}{2}(d_i-1)(d_i-2)$. Hence,
\begin{equation}
\label{eq5} 2\sum_{i=1}^l M_i \leq \sum_{i=1}^l  (d_i-1)(d_i-2)
\leq (\delta-1)(\delta-2).
\end{equation}

The number of points of $\cK$ lying on linear components $r_i$ is
at most $kM_0\leq dM_0$.

For every $\Xi=\Xi_i$, there exists an absolutely irreducible
curve $\Theta$, defined over the algebraic extension $\GF(q^\xi)$
of degree $\xi>1$ of $\GF(q)$ in $\overline{\GF(q)}$, such that
the absolutely irreducible components of $\Xi$ are $\Theta$ and
its conjugates $\Theta_1,\ldots,\Theta_{\xi-1}$. Here, if $\Theta$
has equation $\sum a_{kj}X^kY^j=0$ and $1\leq w \leq \xi-1$, then
$\Theta_w$ is the curve of equation $\sum a_{kj}^{q^w}X^kY^j=0$.
Since $\Xi$, $\Theta$ and the conjugates of $\Theta$ pass through
the same points in $\PG(2,q)$, from B\'ezout's theorem, see
\cite[Lemma 2.24]{H},  $\Xi$ has at most $\theta^2$ points in
$\PG(2,q)$ where $\theta=\deg \Theta$. Note that
$\deg \Xi \geq 2\theta$.

Let $N_i'$ denote the total number of  points (simple or singular) of
$\Xi_i$ lying in $\PG(2,q)$. From the above argument,
\begin{equation}
\label{eq3}
\sum_{i=1}^s N_i'\leq \sum_{i=1}^s \theta_i^2<
(\sum_{i=1}^s\,\theta_i)^2\leq \frac{1}{4}( \sum_{i=1}^s \deg
\Xi_i)^2<\frac{1}{2}(q+2)\sum_{i=1}^s \deg \Xi_i.
\end{equation}
As,
\begin{equation} \label{eq:a}
 qt+\alpha\leq \sum_{i=1}^l (N_i+M_i)+kM_0+\sum_{i=1}^s N'_i,
\end{equation}
from \eqref{eq6}, \eqref{eq5}, \eqref{eq3} and \eqref{eq:a}
it follows that
\[
2(qt+\alpha) \leq
\begin{cases}
2d^2-6d+2+2dM_0+(q+2)d&\mbox{$q\geq 3$ prime\/} \\
  d^3-d^2-6d+2-m_0d(d-3)+2dM_0+(q+2)d&\mbox{otherwise.} \\
\end{cases}
\]
Since  $d\leq 2t-1$, the main assertion follows
by straightforward computation.
\end{proof}
\begin{remark}
As $1\leq m_0\leq M_0\leq d$, the proof of Theorem \ref{teo:main}
shows that Condition \eqref{eq:co1}, and, for $q$ prime, Condition
\eqref{eq2},
may be replaced by the somewhat weaker, but more
manageable, condition $q>8t^3-12^2+4t-2\alpha+2$ (and
$q>16t^2-24t-2\alpha+8$ for $q$ prime).
\end{remark}

\begin{remark}
\label{remarcs} As pointed out in the Introduction, Theorem
\ref{teo:main} is sharp as the bound is attained by some maximal
$(k,n)$-arcs.
\end{remark}

\begin{corollary}\label{lem2}
If $\Gamma$ has a component not defined over $\GF(q)$, then
$$\deg\Gamma \geq 2t+1.$$
\end{corollary}
\begin{proof}
We use the same arguments as in the proof of
Theorem \ref{teo:main}, considering that
$\delta<d$. In particular, we have
\[ 2(qt+\alpha)< d^3-(4+m_0)d^2+(q+5m_0+1)d+(4-4m_0-q)+2dM_0, \]
which for $d\leq 2t$ proves the assertion.
\end{proof}
\begin{remark}
  Corollary \ref{lem2} implies that a plane algebraic curve of
  of degree $2t$ containing $\cK$ is always
  defined over $\GF(q)$.
\end{remark}
\begin{theorem}\label{teo2}
 If the curve $\Gamma$ in Theorem \ref{teo:main} has
 no quadratic component, and
\begin{equation}
  \label{eqq}
  q> \frac{750t^3-1725t^2+10(10M_0+113)t-184-40(\alpha+M_0)}{40},
\end{equation}
 then
 \[ d \geq \frac{5}{2}t. \]
 If $q>5$ is prime, then Condition \eqref{eqq} may be relaxed to
 \begin{equation}
\label{eqqbis}
  q > \frac{125t^2+2(10M_0-105)t-8(\alpha+M_0-9)}{8}
\end{equation}
\end{theorem}
\begin{proof}
 We prove that if $d\leq\frac{5}{2}(t-1)$,
 then \eqref{eqq} (and, for $q$ prime, \eqref{eqqbis})
 does not hold. It is sufficient just
 a change in the proof of Theorem \ref{teo:main}.
 Let $\nu_0^{(i)}=0,$ $\ldots,\nu_4^{(i)}$ be the
 Frobenius orders of $\Delta_i$ with respect to conics. If $q$ is
 a prime greater than $5$, then $\nu_j^{(i)}=j$ for $0\leq j \leq 4$.
 Otherwise, set $\nu^{(i)}=\sum_{j=1}^{4} \nu_j^{(i)}$. Since
 $\nu^{(i)}\leq 2+\varepsilon_3+\varepsilon_4+\varepsilon_5$ and
 $\varepsilon_5\leq 2d_i$, we have that $\nu^{(i)}\leq 6d_i-1$.
 {}From \eqref{svconic},
\begin{equation}
\label{eq8}
5\sum_{i=1}^l N_i\leq \begin{cases}
10\delta(\delta-3)+(q+5)2\delta  & \mbox{$q>5$ prime\/}. \\
   (6\delta-1)\delta(\delta-3)+(q+5)2\delta, &\mbox{otherwise.}\\

\end{cases}
\end{equation}
Using the same argument as in \eqref{eq:a},
\begin{equation}
\label{eq9}
\sum_{i=1}^s N_i'\leq\frac{1}{4}\left(\sum_{i=1}^s\deg\Xi_i\right)^2<
\frac{2}{5}(q+5)\sum_{i=1}^s\deg\Xi_i.
\end{equation}
Using now \eqref{eq5}, \eqref{eq3}, \eqref{eq8} and \eqref{eq9}
we obtain
\begin{equation}
\label{eqpr} (qt+\alpha)\leq\begin{cases}
\frac{5}{2}d^2-\frac{11}{2}d+dM_0+\frac{2}{5}qd+1&\mbox{$q>5$
prime}\\
\frac{6}{5}d^3-\frac{33}{10}d^2+\frac{11}{10}d+dM_0+\frac{2}{5}qd+1 &\mbox{otherwise.}\\
\end{cases}
\end{equation}
Then, \eqref{eqpr}  does not
hold for any $q$. If $q$ is prime, also
\eqref{eqqbis} is not satisfied.
\end{proof}

\section{Algebraic curves passing through
 the  points of a maximal arc}
\label{maximal} Remark \ref{remarcs} motivates the study of plane
algebraic curves passing through all the points of maximal
$(k,n)$--arc in $\PG(2,q)$.

In this section $\cK$ always denotes a maximal $(k,n)$--arc. Recall
that a \emph{$(k,n)$--arc} $\cK$ of a projective plane $\pi$ is a
set of $k$ points, no $n+1$ collinear. Barlotti \cite{ba} proved
that $k\leq (n-1)q+n$, for any $(k,n)$--arc in $\PG(2,q)$; when
equality holds, a $(k,n)$--arc is \emph{maximal}. A purely
combinatorial property characterising  a $(k,n)$--maximal arc $\cK$
is that every line of $\PG(2,q)$ either meets $\cK$ in $n$ points or
is disjoint from it. Trivial examples of maximal arcs  in $\PG(2,q)$
are  the $(q^2+q+1,q+1)$--arc given by all the points of $\PG(2,q)$
and the $(q^2,q)$--arcs consisting of the points of an affine
subplane $\mathrm{AG}(2,q)$ of $\PG(2,q)$. Ball, Blokhuis and
Mazzocca \cite{BBM}, \cite{BB2} have shown that no non--trivial
maximal arc exists in $\PG(2,q)$ for $q$ odd. On the other hand, for
$q$ even, several maximal arcs exists in the Desarguesian plane and
many constructions are known; see \cite{D}, \cite{T74}, \cite{T80},
\cite{P} \cite{M}, \cite{HM}. The arcs arising from these
constructions, with the exception of those of \cite{T74}, see also
\cite{P},  all consist of the union of $n-1$  disjoint
conics together with their common nucleus $N$. In other words, these
arcs are covered by a completely  reducible curve of degree $2n-1$,
whose components are $n-1$ conics and a line through the point $N$.

\begin{remark}
\label{rembis}   From Corollary \ref{lem2}, if $\Gamma$ has a
component defined over $\oGF{q}$ but not over $\GF(q)$, and it
passes through all the points of $\cK$, then its degree $d$ is at
least $2n-1$.
\end{remark}
The following theorem shows that the above hypothesis on the
components of $\Gamma$ can be dropped as far as $q$ is
sufficiently large.

\begin{theorem}\label{teo1.5}
 For any $n$, there exists  $q_0\leq (2n-2)^2$ such that
 if a plane algebraic curve $\Gamma$ defined over $\GF(q)$ with
 $q>q_0$ passes through all the points of a maximal $(k,n)$--arc
 $\cK$ of $\PG(2,q)$ then its degree $d$
is at least $2n-1$. If equality holds then $\Gamma$ has
either one
linear and $n-1$ absolutely irreducible quadratic components or
$n-2$ absolutely irreducible quadratic components and one
cubic component.
\end{theorem}
The proof depends on the the following lemma.
\begin{lemma}\label{teo3}
  Assume that $\Gamma$ is reducible and that the number of
  its components is less than $n-1$.
  Then, the degree
  $d$ of $\Gamma$ satisfies
  \[d \geq \sqrt[4]{q}.\]
\end{lemma}
\begin{proof}
We use the same setup as in the proof of Theorem \ref{teo:main}.
This time the Hasse--Weil bound is used in place of the
St\"ohr-Voloch bound. The number of points of
$\Delta_i$ in $\PG(2,q)$ is $N_i+M_i$; from \eqref{hasse},
\begin{equation}
\label{eq:hw}
 M_i+N_i\leq q+1+(d_i-1)(d_i-2)\sqrt{q}.
\end{equation}
From \eqref{eq3},
\begin{equation}
\label{eq0}
\sum_{i=1}^s N_i'< (q+1)+\sum_{i=1}^s (\deg \Xi_i-1)(\deg\Xi_i-2)\sqrt{q}.
\end{equation}
If $\Gamma$ has $w$ components, then
\begin{equation}\label{deg5}
w(q+1)+\sum_{i=1}^s (\deg\Xi_i-1)(\deg\Xi_i-2)\sqrt{q}+
\sum_{j=1}^l (d_i-1)(d_i-2)\sqrt{q}\geq (n-1)q+n.
\end{equation}
Since $w\leq n-2$, \eqref{deg5} yields
\[(d-1)(d-2)\sqrt{q}\geq q+2;\]
hence, $d\geq \sqrt[4]{q}$.
\end{proof}

\begin{proof}[Proof of Theorem \ref{teo1.5}]
  Suppose $\Gamma$ to have degree $d<2n-1$.
  By Remark \ref{rembis},
  all components of
  $\Gamma$ are defined  over $\GF(q)$.
  If $\Gamma$ is absolutely irreducible, then \eqref{hasse} implies
  that $\Gamma$ contains at most $q+(2n-3)(2n-4)\sqrt{q}+1$ points.
  However, for $q$ large enough,  this number is less than
  $(n-1)q+n$; a contradiction.
  \par
  When $\Gamma$ has more then one component, denote by $t_j$ the
  number of its components of degree $j$. Let $u$ be the maximum
  degree of such components.  Then, $u\leq 2n-3$ and
  \[ d=\sum_{j=1}^{u} jt_j\leq 2n-2. \]
  {}From \eqref{hasse},
  \[ |\Gamma\cap\cK|\leq nt_1+\sum_{j=2}^{u}
  t_j(q+(j-1)(j-2)\sqrt{q}+1)=\left(\sum_{j=2}^{u}t_j\right)
  q+c\sqrt{q}+d,
  \]
  where
  \[ c=\sum_{j=2}^u t_j(j-1)(j-2),\qquad
  d= nt_1+\sum_{j=2}^u t_j. \]
  Both $c$ and $d$ are independent from $q$; therefore,
  \begin{equation}
    \label{eqlm}
   n-1=\lim_{q\to\infty}\frac{|\Gamma\cap\cK|}{q}=\sum_{j=2}^u
  t_j.
  \end{equation}
  Hence,
\[2(n-1)=2\sum_{j=2}^{u}t_j\leq t_1+\sum_{j=2}^{u} jt_j=d. \]
Since $d\leq 2n-2$, by Lemma \ref{teo3}, for $q>(2n-2)^4$
the curve $\Gamma$ should have at least $n-1$ components.
This would imply that either $u=2$, $t_1=1$, $t_2=(n-1)$
or $u=3$, $t_1=0$, $t_2=(n-2)$, $t_3=1$.
In particular, in both cases  $d=2n-1$.
\end{proof}
\begin{remark}
\label{remwp} As mentioned in the Introduction, case
$d=2n-1$ in Theorem \ref{teo1.5} occurs  when $\cK$ is a Denniston
maximal arc \cite{D} (or one of the maximal arcs constructed by
Mathon and others, \cite{T74,T80,P,M,HM}). This result may not
extend to any of the other known maximal arcs; they are the Thas
maximal $(q^3-q^2+q,q)$--arcs in $\PG(2,q^2)$ arising from the
Suzuki--Tits ovoid of $\PG(3,q)$; see \cite{T74}. In fact, $22$ is
the minimum degree of a plane curve which passes through all
points of a Thas' maximal $(456,8)$--arc in $\PG(2,64)$; see
\cite{webp}.
\end{remark}

\section{Maximal arcs of degree $4$}
\label{deg4}
 From Theorem \ref{teo1.5}, for $q>6^4$,
 a lower bound on the degree
 of an algebraic curve $\Gamma$ passing through all the points of a
 maximal arc $\cK$ of degree $4$ is $7$. Our aim  is to prove in
 this case the following result.


\begin{theorem}\label{teo5}
 Let $\cK$ be a maximal arc of degree $4$ and suppose there
 exists an algebraic plane curve $\Gamma$ containing all the
 points of $\cK$. If $\deg\Gamma=7$, then
 $\Gamma$ consists of three
 disjoint conics, all with the same nucleus $N$, and a line through
 $N$.
\end{theorem}
\begin{proof}
From Theorem \ref{teo1.5}, the curve $\Gamma$  splits either into
one irreducible cubic  and two irreducible conics or into three
irreducible conics and one line $r$. These two cases are
investigated separately.

Let $C$ denote any of the above conics of nucleus $N$.
We show that every
point of $C'=C\cup\{N\}$ lying in $\PG(2,q)$ is contained in $\cK$:
in fact,
if there were a point $P\in C'\setminus \cK$, then there would be at
least $\frac{q}{4}$ lines through $P$ external to $\cK$. All these
lines would meet $C'$ in $\frac{q}{4}$
distinct points, which, in turn, would
not be on $\cK$.
Hence,  $\Gamma$ would have  less than $3q+4$ points
on the arc $\cK$, a contradiction.

Now assume that $\Gamma$ splits into a  cubic $\cD$ and two conics
$C_i$, with $i=1,2$. Denote by $N_i$ the nucleus of $C_i$ and set
$\cX=C_1 \cup C_2 \cup\{N_1, N_2\}$. Since $|\cX|\leq 2q+4$, there
exists a point $P\in \cK \setminus \cX$. Obviously, $P\in\cD$. Every
line through $P$ meets $\cK$ in four points; thus, there is no line $\ell$
through $P$ meeting  both $C_1$ and $C_2$ in $2$ points;
otherwise, $|\ell\cap\cX\cap\cK|=4$ and $|\ell\cap\cK|\geq 5$,
a contradiction.
Hence, there are at least
$q-1$ lines through $P$ meeting $\cD$ in another point $P'$.
There are at most $5$ bisecants to the irreducible cubic
curve $\cD$ through any given point $P\in\cD$, namely the
tangent in $P$ to $\cD$ and, possibly, four other tangents
in different points to $\cD$ passing through $P$.
Hence, there are
$q-6$ lines through $P$ meeting $\cD$ in three points.
If this were the
case, $\cD$ would consist of at least $2(q-6)+6$ points, which is
impossible.
Therefore, we may assume that $\Gamma$ splits into three conics, say
$C_1$, $C_2$, $C_3$, with nuclei $N_1$, $N_2$, $N_3$, and a line
$r$.

Recall that, as seen above,
 the nuclei of all the conics belong to $\cK$.
Now we  show that at least one nucleus, say $N_1$, lies on the line
$r$. Since $\Gamma$ is a curve containing $\cK$ of minimum degree
with respect to this property, there is at least  a point $P$ on $
r\cap \cK$ not on $C_1$. Each line through $P$ is a $4$--secant to
$\cK$ hence, it meets $\cX$  in an odd number of points. If $P\neq
N_1$, then the number of lines through $P$ meeting $\cX$ in an odd
number of points is at most $15$, which is less than $q+1$ for
$q\geq 2^4$, a contradiction.
Actually, all the nuclei $N_i$  lie on $r$. In fact, suppose that
$N_j \notin r$ for $j\in\{2,3\}$. Then, $N_j \in C_s$, with $s\neq
1,j$ and the line ${N_1N_j}$ joining $N_1$ and $N_j$ is
tangent to $C_1$ and $C_j$. Consequently, ${N_1N_j}$ meets
$C_s$ in another point different from $N_j$ that is, it is a
$5$--secant to $\cK$, again a contradiction.
\par
We are left with three cases, namely:
\begin{enumerate}[(1)]
\item
\label{i:1}
 $N_1\neq N_2 \neq N_3$, $C_i \cap C_j=\{A\}$ for any $i\neq j$ and
  $A \in N_1N_2$.
\begin{figure}
\[
\psfrag{A}{\footnotesize $A$} \psfrag{n1}{\footnotesize $N_1$}
\psfrag{n2}{\footnotesize $N_2$} \psfrag{n3}{\footnotesize $N_3$}
\psfrag{c1}{\footnotesize $C_1$} \psfrag{c2}{\footnotesize
$C_2$} \psfrag{c3}{\footnotesize $C_3$}
 \includegraphics[width=6cm]{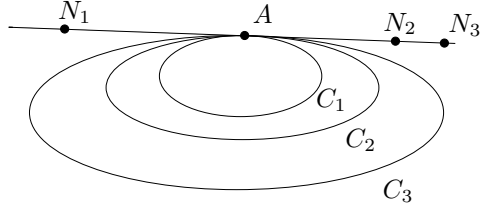} \]
\caption{Case \ref{i:1} in Theorem \ref{teo5}}
\end{figure}
\item
\label{i:2}
 $N_i= N_j, \ N_j\neq N_s$, $N_i\in C_s$,
  $C_i \cap C_j = \emptyset $,
  with  $i,j,s\in \{1,2,3\}$.
\begin{figure}
\[
\psfrag{ci}{\footnotesize $C_i$} \psfrag{cj}{\footnotesize
$C_j$} \psfrag{cs}{\footnotesize $C_s$}
\psfrag{ninj}{\footnotesize $N_i=N_j$} \psfrag{ns}{\footnotesize
$N_s$}
 \includegraphics[width=8cm]{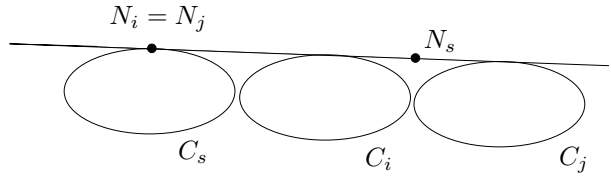} \]
\caption{Case \ref{i:2} in Theorem \ref{teo5}}
\end{figure}
\item \label{i:3}
 $N_1=N_2=N_3$ and $C_i \cap C_j = \emptyset$ for $i\neq j$.
\begin{figure}
\[
\psfrag{ci}{\footnotesize $C_1$} \psfrag{cj}{\footnotesize
$C_2$} \psfrag{cs}{\footnotesize $C_3$}
\psfrag{nj}{\footnotesize $N_1=N_2=N_3$}
 \includegraphics[width=10cm]{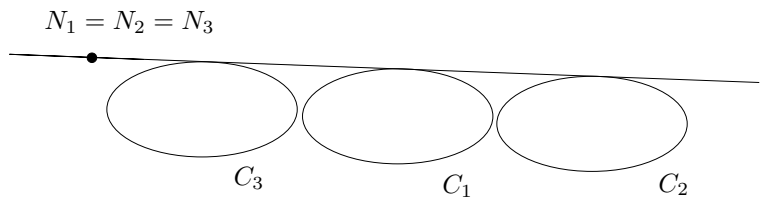} \]
\caption{Case \ref{i:3} in Theorem \ref{teo5}}
\end{figure}
\end{enumerate}
We are going to show that cases (\ref{i:1}) and (\ref{i:2}) do not
actually occur.

\begin{lemma}\label{lemma2}
  Let $C_1$, $C_2$ be two  conics with a common point $A$ but different
  nuclei $N_1$, $N_2$. If $A\in {N_1N_2}$  then  there is a line $s$ with $A\in s$
  and  $s \neq {N_1N_2}$ such that for any point on $s\setminus
  (C_1\cup C_2)$ there passes a line $\ell$  with $$|\ell \cap (C_1 \cup
  C_2 \cup\{N_1, N_2\})|\geq 3.$$
\end{lemma}
\begin{proof}
Let $(X,Y,Z)$ denote homogeneous coordinates of  points of the plane
$\PG(2,q)$. Choose a reference system such that $A=O=(0,0,1)$ and the
line joining $N_1$ and $N_2$ is the $X$--axis. We may suppose $C_i$
to have equation
\[\alpha_iX^2+XY+\beta_iY^2+\lambda_iYZ=0, \]
where $\alpha_i, \beta_i,\lambda_i \in\GF(q)$ and $i=1,2$.

Since both $C_i$ are non--degenerate conics, we have $\alpha_i \neq
0$ and $ \lambda_i \neq 0$. Furthermore, $\lambda_1\neq \lambda_2$
as the nuclei $N_1$ and $N_2$ are distinct.

Denote by $\tr(x)$ the trace of $\GF(q^2)$ over $\GF(q)$;
namely $\tr(x)=x+x^q$.
If $\tr((\alpha_1\lambda_2+\alpha_2\lambda_1)(\beta_1\lambda_2+
\beta_2\lambda_1))=0$, then the two conics $C_1$ and  $C_2$ have
more than one point in common, which is impossible.  Hence,
 $\tr(\alpha_1\lambda_2+\alpha_2\lambda_1)(\beta_1\lambda_2+
\beta_2\lambda_1))=1$;  in particular,
$\beta_1\lambda_2\neq \beta_2\lambda_1$.

A generic point $P_m^i$ of $C_i\setminus \{A\}$  has
homogeneous coordinates
\[P_m^i=\Big(\frac{\lambda_i m}{\alpha_i m^2+m+\beta_i},\frac{\lambda_i
}{\alpha_i m^2+m+\beta_i},1\Big),\] with $m\in\GF(q) \setminus
\{0\}$. Consider now a point  $P_{\varepsilon}=(0,\epsilon, 1)$ on
the $Y$--axis,  with $\varepsilon \in\GF(q) \setminus \{ 0,
\frac{\lambda_1}{\beta_1}, \frac{\lambda_2}{\beta_2} \}$. The
points $P_m^1$, $P_t^2$ and $P_{\varepsilon}$ are collinear if and
only if
\[ \left|\begin{array}{ccc}
    {{\displaystyle\frac{\lambda_1 m}{\alpha_1 m^2+m+\beta_1}}} &
{\displaystyle\frac{\lambda_1}{\alpha_1 m^2+m+\beta_1 }}& 1 \\
{\displaystyle\frac{\lambda_2t}{\alpha_2 t^2+t+\beta_2}}
 & {\displaystyle\frac{\lambda_2}{\alpha_2 t^2+t+\beta_2} }& 1
\\0 & \varepsilon & 1
\end{array}
\right|=0, \] that is
\begin{equation}\label{cubic}
  \alpha_1\lambda_2\varepsilon m^2t+\alpha_2\lambda_1\epsilon m
  t^2+ \varepsilon(\lambda_2+\lambda_1)mt+
  (\varepsilon \lambda_1 \beta_2+\lambda_1\lambda_2)m+(\epsilon
  \lambda_2\beta_1+\lambda_1\lambda_2)t=0.
\end{equation}
Equation \eqref{cubic} may be regarded as the affine equation of a
cubic curve $\cD$ in the indeterminate $m$ and $t$. Observe that
$(0,0,1) \in \cD$. The only points at infinity of $\cD$ are
$Y_{\infty}=(0,1,0)$, $X_{\infty}=(0,1,0)$ and
$B=(\frac{\lambda_1\alpha_2}{\lambda_2\alpha_1},1,0)$. Therefore,
$\cD$ does not split into three conjugate complex lines. Thus, by
\cite[Theorem 11.34]{H} and \cite[Theorem 11.46]{H}, there is at
least one affine point $T=(\overline{m},\overline{t},1)$ on $\cD$
different from $(0,0,1)$.

Since $\tr((\alpha_1\lambda_2+\alpha_2\lambda_1)(\beta_1\lambda_2+
\beta_2\lambda_1))=1$, the line $m=t$ is a $1$--secant to $\cD$ in
$(0,0,1)$; hence,  the point $T$ is not on this line.

 This implies that, for any given
$\varepsilon \in\GF(q) \setminus \{0\}$, there exist at least two
distinct values $\overline{m}, \ \overline{t} \in\GF(q) \setminus
\{0\}$ satisfying \eqref{cubic}. Hence, $P_{\overline{m}}^1$,
$P_{\overline{t}}^2$ and $P_{\varepsilon}$ are collinear and the
line $P_{\overline{m}}^1P_{\overline{t}}^2$ meets $C_1 \cup
  C_2 \cup\{N_1, N_2\}$ in at least
three points.
\end{proof}

{}From Lemma $\ref{lemma2}$, in case (\ref{i:1}) the set $C_1 \cup
C_2 \cup\{N_1, N_2\}$   cannot be completed
to a maximal arc just by adding a third conic $\cC_3$, together
with its nucleus $N_3$, since, in this case, there would be at
least a $5$--secant to $C_1\cup C_2\cup C_3\cup\{N_1,N_2,N_3\}$.
Hence, case (\ref{i:1}) is ruled out.

\begin{lemma}
  \label{lemma1}
  Given any two disjoint conics $C_1$, $C_2$  with the same nucleus
  $N$, there is a unique  degree--$4$ maximal arc containing
  $\cX=C_1\cup C_2$.
\end{lemma}
\begin{proof}
There is a line $r$ in $\PG(2,q)$ external to $\cX$, since,
otherwise, $\cX$ would be a $2$--blocking set with less than
$2q+\sqrt{2q}+1$ points, which is a contradiction; see \cite{BB}.

Choose a reference system such that $N=O=(0,0,1)$ and $r$
is the line at infinity $Z=0$.  The conics $C_i$, for $i=1,2$,
have equation:
\begin{equation}
\alpha_iX^2+XY+\beta_iY^2+\lambda_iZ^2=0,
\end{equation}
where $\alpha_i, \beta_i, \lambda_i \in \GF(q)$ and
$\tr(\alpha_i\beta_i)=1$.

Since both $C_i$ are  non--degenerate,  $\lambda_i\neq 0$.
We first show that,
$\lambda_1 \neq \lambda_2$, as  $C_1$ and $C_2$ are disjoint.
We argue by contradiction. If it were
 $\lambda_1=\lambda_2$, then we could assume
 $\alpha_1 \neq \alpha_2$; in fact, if $\alpha_1=\alpha_2$
 and $\lambda_1=\lambda_2$,  the linear
 system generated by $C_1$ and $C_2$ would contain the
 line $Y=0$; thus their intersection would not be
 empty.
Let now
$$\gamma=
\sqrt{\frac{\beta_1-\beta_2}{\alpha_1-\alpha_2}};$$
hence,
\[\alpha_1 \gamma^2+\gamma+\beta_1=
\alpha_2\gamma^2+\gamma+\beta_2,\]
and the  line $X=\gamma Y$
would meet  the two conics in the  same  point
\[ P=\Big(\gamma\big({\textstyle
  \frac{\lambda_1}{\alpha_1 \gamma^2+\gamma+\beta_1}}\big)^{\frac{1}{2}},
\big({\textstyle
  \frac{\lambda_1}{\alpha_1 \gamma^2+\gamma+\beta_1}}\big)^{\frac{1}{2}},1
\Big), \]
contradicting $C_1 \cap C_2= \emptyset.$

We also see that
$\alpha_1\lambda_2\neq \alpha_2\lambda_1$  and $\beta_1\lambda_2\neq
\beta_2\lambda_1$, since,  otherwise, the points
\[
P=
({\textstyle \sqrt{\frac{\lambda_1}{\alpha_1}}},0,1)=
({\textstyle \sqrt{\frac{\lambda_2}{\alpha_2}}},0,1),
\quad
Q=
(0,{\textstyle\sqrt{\frac{\lambda_1}{\beta_1}}},1)=
(0,{\textstyle\sqrt{\frac{\lambda_2}{\beta_2}}},1)\]
 would lie on both conics.

Let now $C_3$ be the conic with equation
\[\frac{\alpha_1\lambda_2+
\alpha_2\lambda_1}{\lambda_1+\lambda_2}X^2+XY+\frac{\beta_1\lambda_2+
\beta_2\lambda_1}{\lambda_1+\lambda_2}Y^2+(\lambda_1+\lambda_2)Z^2=0.\]
Set
\[\nu=\frac{(\alpha_1\lambda_2+\alpha_2\lambda_1)(\beta_1\lambda_2+%
\beta_2\lambda_1)}{\lambda_1^2+\lambda_2^2}. \]
The collineation H
of $\PG(2,q)$ given by the matrix
$$\left( \begin{array}{ccc} a^{-1} & 0 & 0 \\ 0 & a & 0
\\ b & c & 1
\end{array} \right), $$
where $a=\sqrt{\frac{\alpha_1\lambda_2+
\alpha_2\lambda_1}{\lambda_1+\lambda_2}}, \
b=\sqrt{\frac{1+\alpha_1/a^2}{\lambda_1}}$ and
$c=a\sqrt{\frac{\beta_1+\beta_2}{\lambda_1+\lambda_2}},$ maps the
conics $C_i$,  $i=1,2$, to
\[\overline{C}_2^i: X^2+XY+\nu Y^2+\lambda_iZ^2=0,\]
and $C_3$ to
 \[\overline{C}_2^3: X^2+XY+\nu Y^2+(\lambda_1+\lambda_2)Z^2=0.\]
If it were $\tr(\nu)=0$, then $\overline{C}_2^1$ and
$\overline{C}_2^2$ would share some point in common on the line at
infinity. Hence, $\tr(\nu)=1$.  In particular, $\overline{C}_2^1$,
$\overline{C}_2^2$ and  $\overline{C}_2^3$  together with their
common nucleus $O$, form a degree--$4$ maximal arc $\cK$ of Denniston
type; see \cite{D}, \cite{al} and \cite[Theorem 2.5]{M}.

It remains to show the uniqueness of $\cK$.  We first observe that
$\cX$ is a $(0,2,4)$--set with respect to lines of the plane. No
point lying on a $4$--secant to $\cX$ can  be  added to  $\cX$ to
get a degree--$4$ maximal arc.

Take $P \notin \cX$ and  denote by
  $u_i$ with $i=0,2,4$,  the number of $i$--secants to $\cX$ through
  $P$,
that is the number of lines meeting $\cX$ in $i$ points.
The lines through $P$  which are external to $\cX$ are  also
external $\cK$ and the converse also holds.   Therefore, when
$P\not\in\cK$, we
  have $u_0=\frac{1}{4}q$. From
 \[ u_0+u_2+u_4=q+1,\]
 \[2u_2+4u_4=2q+2,\]

 also $u_4=\frac{1}{4}q$.
 Hence, no point $P\not\in\cK$ may be added to
 $\cX$  to obtain a maximal arc of degree $4$.
\end{proof}
Finally, Lemma \ref{lemma1} shows that case (\ref{i:2}) does not
occur.
\end{proof}

\penalty-10
\vskip.5cm\noindent {\em Authors' addresses}:\\
\penalty10000
\noindent\vskip.2cm
\penalty10000
\begin{minipage}{6cm}
  Angela AGUGLIA \\
  Dipartimento di Matematica \\
  Politecnico di Bari \\
  Via Orabona 4 \\
  70125 Bari  (Italy) \\
  E--mail: {\tt a.aguglia@poliba.it} \\
\end{minipage}
\hfill
\begin{minipage}{6cm}
  Luca GIUZZI \\
  Dipartimento di Matematica \\
  Politecnico di Bari \\
  Via Orabona 4 \\
  70125 Bari  (Italy) \\
  E--mail: {\tt l.giuzzi@poliba.it} \\
\end{minipage}
\vskip.1mm
\begin{minipage}{6cm}
  G\'abor KORCHM\'AROS\\
  Dipartimento di Matematica\\
  Universit\`a della Basilicata\\
  Contrada Macchia Romana\\
  85100 Potenza (Italy).\\
  E--mail: {\tt korchmaros@unibas.it}
\end{minipage}

\end{document}